\input amstex
\input amsppt.sty
\magnification=\magstep1
\vsize=22.2truecm
\baselineskip=16truept
\NoBlackBoxes
\pageno=1
\nologo
\def\Z{\Bbb Z}
\def\N{\Bbb N}

\def\l{\left}
\def\r{\right}
\def\bg{\bigg}
\def\({\bg(}
\def\[{\bg\lfloor}
\def\){\bg)}
\def\]{\bg\rfloor}
\def\t{\text}
\def\f{\frac}

\def\p{\ (\roman{mod}\ p)}

\def\bi{\binom}
\def\eq{\equiv}
\def\cs{\cdots}
\def\ls{\leqslant}
\def\gs{\geqslant}
\def\mo{\roman{mod}}

\def\ve{\varepsilon}
\def\al{\alpha}
\def\da{\delta}

\def\st#1#2{\thickfrac\thickness0{#1}{#2}}
\def\nh{\noalign{\hrule}}
\def\hh{height4pt}

\def\om{&\omit &}
\def\M#1#2{\thickfracwithdelims[]\thickness0{#1}{#2}_m}

\def\Mm#1#2#3{\thickfracwithdelims[]\thickness0{#1}{#2}_{#3}}

\def\m#1#2{\thickfracwithdelims\{\}\thickness0{#1}{#2}_m}
\def\mm#1#2#3{\thickfracwithdelims\{\}\thickness0{#1}{#2}_{#3}}
\def\bm#1#2#3{\thickfracwithdelims ()\thickness0{#1}{#2}_{#3}}
\def\bmb#1#2#3{\thickfracwithdelims ()\thickness0{#1}{#2}^*_{#3}}
\def\Proof{\noindent{\it Proof}}

\def\Remark{\medskip\noindent{\it  Remark}}

\def\Ack{\medskip\noindent {\bf Acknowledgment}}
\hbox {Discrete Math. 308(2008), no.\,18, 4231--4245.}
\bigskip
\topmatter
\title On Sums of Binomial Coefficients and their Applications\endtitle
\author Zhi-Wei Sun\endauthor
\affil Department of Mathematics, Nanjing University
\\ Nanjing 210093, People's Republic of China
\\zwsun\@nju.edu.cn
\\{\tt http://math.nju.edu.cn/\~{}zwsun}
\endaffil
\abstract In this paper we study recurrences concerning the
combinatorial sum $\M nr=\sum_{k\eq r\ (\mo\ m)}\bi nk$ and the
alternate sum $\sum_{k\eq r\ (\mo\ m)}(-1)^{(k-r)/m}\bi nk$,
where $m>0$, $n\gs0$ and $r$ are integers. For example, we show
that if $n\gs m-1$ then
$$\sum_{i=0}^{\lfloor(m-1)/2\rfloor}(-1)^i\bi{m-1-i}i\M{n-2i}{r-i}=2^{n-m+1}.$$
We also apply such results to investigate Bernoulli and Euler
polynomials. Our approach depends heavily on an identity
established by the author [{\it Integers} {\bf 2}(2002)].

\endabstract
\thanks 2000 {\it Mathematics Subject Classifications}. Primary 11B65;
Secondary 05A19, 11B37, 11B68.
\newline\indent The initial version of this paper was posted as
{\tt arXiv:math.NT/0404385} on April 21, 2004.
\newline\indent The author was supported by
the National Science Fund for Distinguished Young Scholars in
China (Grant No. 10425103).
\endthanks
\endtopmatter
\document

\heading{1. Introduction and Main Results}\endheading

As usual, we let
$$\bi x0=1\ \t{and}\ \bi xn=\f{x(x-1)\cdots(x-n+1)}{n!}\ \t{for}\ n\in\Z^+=\{1,2,3,\ldots\}.$$
Following [Su2], for $m\in\Z^+$, $n\in\N=\{0,1,2,\ldots\}$ and $r\in\Z$
we set
$$\M{n}r=\sum^n\Sb k=0\\k\eq r\ (\mo\ m)\endSb \bi nk
\ \ \t{and}\ \
\m nr=\sum^n\Sb k=0\\k\eq r\ (\mo\ m)\endSb (-1)^{\f{k-r}m}\bi nk.\tag1.1$$
As $\bi {n+1}{k}=\bi nk+\bi n{k-1}$
for any $k\in\Z^+$, we have the following useful recursions:
$$\M{n+1}{r}=\M nr+\M n{r-1}\ \ \t{and}\ \ \m{n+1}{r}=\m nr+\m n{r-1}.\tag1.2$$

Let $m,n\in\Z^+$ and $r\in\Z$. The study of the sum $\M nr$ dates back to 1876 when
C. Hermite  showed that if $n$ is odd and $p$ is an odd prime then
$\Mm n0{p-1}\eq 1\ (\mo\ p)$ (cf. L. E. Dickson [D, p.\,271]). In 1899 J. W. L.
Glaisher obtained the
following generalization of Hermite's result:
$$\Mm{n+p-1}r{p-1}\eq\Mm nr{p-1}\ (\mo\ p)\ \quad \t{for any prime}\ p.$$
(See, e.g., [Gr, (1.11)].) If $p$ is a prime with $p\eq 1\ (\mo\
m)$, then $\M{p}{r}\eq\M{1}{r}\ (\mo\ p)$ since $p$ divides any of
$\bi p1,\ldots,\bi p{p-1}$, thus $\M{n+p-1}{r}\eq\M{n}r\ (\mo\ p)$
by (1.2) and induction. This explains Glaisher's result in a
simple way. (Recently the author and R. Tauraso [ST] obtained a further extension of Glaisher's congruence.)
In the modern investigations made
by Z. H. Sun and the author (cf. [SS], [S], [Su1] and [Su2]), $\M
nr$ was expressed in terms of linear recurrences and then applied
to produce congruences for primes. The sum $\M nr$ also appeared
in C. Helou's study of Terjanian's conjecture concerning Hilbert's
residue symbol and cyclotomic units (cf. [H, Prop. 2 and Lemma 3]).

 Now we state two theorems on the sums in (1.1) and give two corollaries.
 The proofs of them depend heavily on an identity
 established by the author in [Su3], and will be presented in Section 2.

\proclaim{Theorem 1.1} Let $m$ be a positive integer. Then, for any
integers $k$ and $n\gs2\lfloor(m-1)/2\rfloor$, we have
$$\sum_{i=0}^{\lfloor(m-1)/2\rfloor}(-1)^i\bi{m-1-i}i\M{n-2i}{k-i}
=2^{n-m+1}+\da_{m-2,\,n}\f{(-1)^k}2,\tag1.3$$ where the Kronecker
symbol $\da_{l,n}$ is $1$ or $0$ according to whether $l=n$ or not.
\endproclaim

\proclaim{Corollary 1.1} Let $k\in\Z$ and $m\in\Z^+$. For $n\in\N$ set
$$u_n=\M n{\lfloor (k+n)/2\rfloor}\ \t{and}\
v_n=mu_n-2^n-\da_{n,0}\da_{(-1)^m,1}(-1)^{\lfloor k/2\rfloor},\tag1.4$$
where $\lfloor\al\rfloor$ denotes the integral part of a real number $\al$.
Then we have
$$\sum_{i=0}^{\lfloor(m-1)/2\rfloor}(-1)^i\bi{m-1-i}iu_{n-2i}
=2^{n-m+1}-\da_{m-2,\,n}\f{(-1)^{\lfloor
(k+m)/2\rfloor}}2\tag1.5$$ for every integer
$n\gs2\l\lfloor(m-1)/2\r\rfloor$.
Also, $(v_n)_{n\in\N}$ is a
linear recurrence sequence, satisfying the recurrence:
$$\sum_{i=0}^{\lfloor(m-1)/2\rfloor}(-1)^i\bi{m-1-i}iv_{n-2i}=0
\quad\t{for all}\ n\gs 2\l\lfloor\f{m-1}2\r\rfloor.\tag1.6$$
\endproclaim

\Remark\ 1.1. (a) In fact, the author first proved (1.6) in the case $2\nmid m$ on August 1, 1988, motivated by a conjecture
of Z. H. Sun; after reading the author's initial proof Z. H. Sun [S]
noted that the equality in (1.6) also holds if
$2\mid m$ and $n\gs m-1$.
(b) In light of the first equality in (1.2), on August 11, 1988
the author obtained the following result by induction:
Let $m,n\in\Z^+$ and $m>2$. If $n\gs m-1$ then
 $$\M n{\lfloor\f{n+1}2\rfloor}>\M n{\lfloor\f{n+1}2\rfloor+1}>\cdots>\M n{\lfloor\f{n+m}2\rfloor},$$
 otherwise
 $$\M n{\lfloor\f{n+1}2\rfloor}>\cdots>\M nn>\M n{n+1}=\cdots=\M n{\lfloor\f{n+m}2\rfloor}=0.$$
Therefore
$$\M n{\lfloor \f n2\rfloor}>\f{2^n}m>\M n{\lfloor\f{m+n}2\rfloor}.$$
\medskip

\proclaim{Theorem 1.2} Let $k\in\Z$ and $m\in\Z^+$. Then
$$\sum^{\lfloor(m+1)/2\rfloor }_{i=0}(-1)^{i}c_m(i)\M{n-2i}{k-i}
=2(-1)^{k}\da_{m,\,n}\tag1.7$$
for each integer $n\gs2\lfloor(m+1)/2\rfloor $,
and
$$\sum^{\lfloor m/2\rfloor }_{i=0}(-1)^{i}d_m(i)\m{n-2i}{k-i}
=(-1)^{k}\da_{m-1,\,n}\tag1.8$$
for any integer $n\gs2\lfloor m/2\rfloor $, where $c_1(1)=4$, and
$$c_m(i)=\f{m^2+m-2i}{(m-i)(m+1-i)}\bi{m+1-i}i\in\Z\ \t{and}\
d_m(i)=\f{m}{m-i}\bi{m-i}i\in\Z$$
for every $i=0,\ldots,m-1$.
\endproclaim

\Remark\ 1.2.
Let $p$ be an odd prime. It is easy to check that
$$(-1)^{i-1}c_{p-1}(i)\eq(-1)^id_{p-1}(i)\eq C_i\ (\mo\ p)
\ \ \t{for}\ i=1,2,\cs,\f{p-1}2,$$
where $C_i=\bi {2i}i/(i+1)=\bi{2i}i-\bi{2i}{i+1}$ is the $i$-th Catalan number.
\medskip

\proclaim{Corollary 1.2 {\rm (A. Fleck, 1913)}} Let $n\in\Z^+$ and $r\in\Z$. If $p$ is a prime, then
$$\sum\Sb 0\ls k\ls n\\p\mid k-r\endSb(-1)^k\bi nk\eq0\ \l(\mo\ p^{\lfloor(n-1)/(p-1)\rfloor}\r).\tag1.9$$
\endproclaim

By [Su2, Remark 2.1], for $k\in\Z$, $l\in\N$, $m\in\Z^+$ and $\ve\in\{1,-1\}$, we have
$$\sum_{\gamma^m=\ve}\gamma^k(2-\gamma-\gamma^{-1})^l
=(-1)^km\times\cases\M{2l}{k+l}&\t{if}\ \ve=(-1)^m,
\\\m{2l}{k+l}&\t{otherwise}.\endcases$$
So Theorem 1.2 is closely related to
the following materials on Bernoulli and Euler polynomials.

Let $m,n\in\Z^+$, $q\in\Z$ and $(q,m)=1$, where $(q,m)$ is the
greatest common divisor of $q$ and $m$. If $\gamma^m=1$ and
$\gamma^{(q,m)}=\gamma\not=1$, then $\gamma^q\not=1$ and hence
$2-\gamma^q-\gamma^{-q}\not=0$. We define a linear recurrence
$(U^{(q)}_l(m,n))_{l\in\N}$ of order $\lfloor m/2\rfloor$ by
$$U^{(q)}_l(m,n)=\f1{2m}\sum\Sb \gamma^m=1\\\gamma\not=1\endSb
\f{2-\gamma^{qn}-\gamma^{-qn}}{2-\gamma^q-\gamma^{-q}}(2-\gamma-\gamma^{-1})^l.\tag1.10$$
Note that $U_l^{(-q)}(m,n)=U_l^{(q)}(m,n)$ and
$$mU_l^{(q)}(m,n)=(1-(-1)^{(m-1)n})2^{2l-2}+\sum\Sb d\mid m\\d>2\endSb u_l^{(q)}(d,n),$$
where
$$\align u_l^{(q)}(d,n)=&\sum\Sb 0<c<d/2\\(c,d)=1\endSb\f{2-e^{2\pi i\f cdqn}-e^{-2\pi i\f cdqn}}
{2-e^{2\pi i\f cdq}-e^{-2\pi i\f cdq}}\l(2-e^{2\pi i\f cd}-e^{-2\pi i\f cd}\r)^l
\\=&\sum\Sb 0<c<d/2\\(c,d)=1\endSb
\(\f{\sin(\pi nqc/d)}{\sin(\pi qc/d)}\)^2\l(4\sin^2\f{\pi c}d\r)^l.
\endalign$$
Obviously $u_l^{(q)}(m,n)=mU_l^{(q)}(m,n)$ if $m$ is an odd prime.
Later we will see that $U_0^{(q)}(m,n)=n(m-n)/(2m)$ if $1\ls n\ls m$, and
$U_l^{(q)}(m,n)\in\Z$ if $l>0$.
When $(q,2m)=1$, for $l\in\N$ we also define
$$V^{(q)}_l(m,n)=\f1{2m}\sum_{\gamma^m=-1}
\f{2-\gamma^{qn}-\gamma^{-qn}}{2-\gamma^q-\gamma^{-q}}(2-\gamma-\gamma^{-1})^l;\tag1.11$$
clearly
$V^{(\pm q)}_l(m,n)=2U^{(q)}_l(2m,n)-U^{(q)}_l(m,n)$
since $\gamma^m=-1$ if and only if $\gamma^{2m}=1$ but $\gamma^m\not=1$.

 Let $p$ be an odd prime,
 and let $m,n>0$ be integers with $p\nmid m$ and $m\nmid n$.
A. Granville and the author [GS, pp.\,126--129]
proved the following surprising result for Bernoulli
polynomials: If $p\eq \pm q\ (\mo\ m)$ where $q\in\Z$, then
$$B_{p-1}\l(\l\{\f {pn}m\r\}\r)-B_{p-1}\eq\f m{2p}\l(U_p^{(q)}(m,n)-1\r) \ (\mo\ p)\tag1.12$$
where we use $\{\al\}$ to denote the fractional part of a real
number $\al$. (The reader may consult [Su4] for other congruences
concerning Bernoulli polynomials.) With the help of Theorem 1.2,
we can write the recurrent coefficients of the sequence
$(U_l^{(q)}(m,n))_{l\in\N}$ in a simple closed form.

 \proclaim{Theorem 1.3} Let $m,n\in\Z^+$, $q\in\Z$ and $(q,m)=1$.
Then we have the recursions:
$$U^{(q)}_l(m,n)=\sum_{0<i\ls\lfloor  m/2\rfloor }(-1)^{i-1}a_m(i)U^{(q)}_{l-i}(m,n)
\quad\t{for}\ l\gs\l\lfloor  \f m2\r\rfloor, \tag1.13$$
and
$$V^{(q)}_l(m,n)=\sum_{j=1}^{\lfloor(m+1)/2\rfloor }(-1)^{j-1}b_m(j)
V^{(q)}_{l-j}(m,n)
\quad\t{for}\ l\gs\l\lfloor  \f{m+1}2\r\rfloor \tag1.14$$
provided $(q,2m)=1$,
where the integers $a_m(i)$ and $b_m(j)$ are given by
$$a_m(i)=\cases c_m(i)&\t{if}\ 2\mid m,\\
d_m(i)&\t{if}\ 2\nmid m,\endcases\ \
\t{and}\ \ b_m(j)=\cases d_m(j)&\t{if}\ 2\mid m,\\
c_m(j)&\t{if}\ 2\nmid m.\endcases
\tag1.15$$
If $m$ does not divide $n$, and $p$ is an odd prime with $p\eq \pm q\ (\mo\ 2m)$, then
$$(-1)^{\lfloor pn/m\rfloor }E_{p-2}\l(\l\{\f{pn}m\r\}\r)+\f{2^p-2}p
\eq\f mp\l(V^{(q)}_p(m,n)-1\r)
\ (\mo\ p),\tag1.16$$
where the Euler polynomials $E_k(x)\ (k=0,1,\ldots)$ are given by
$$\f{2e^{xz}}{e^z+1}=\sum_{k=0}^{\infty}E_k(x)\f{z^k}{k!}.$$
\endproclaim

In Section 3 we will first deduce Theorem 1.3 from Theorem 1.2, and then
give another proof of (1.13) via Chebyshev polynomials.
Section 4 is an appendix containing the explicit values of $a_m(i)$ and $b_m(j)$
for $m=2,3,\ldots,12$.

\heading{2. Proofs of Theorems 1.1 and 1.2}\endheading
\proclaim{Lemma 2.1} Let $l$ be any nonnegative integer. Then
$$\sum_{j=0}^l(-1)^{l-j}\bi{x+y+j}{l-j}\bi{y+2j}j=\sum_{j=0}^l\bi{l-x}j.\tag2.1$$
\endproclaim
\Proof. Since both sides of (2.1) are polynomials in $x$ and $y$,
it suffices to show (2.1) for all $x\in\{l,l+1,\ldots\}$ and
$y\in\{0,2,4,\ldots\}$.

Let $x=l+n$ and $y=2k$ where $n,k\in\N$. Set $m=k+l$. Then
$$\align&\sum_{j=0}^l(-1)^{l-j}\bi{x+y+j}{l-j}\bi{y+2j}j
\\=&\sum_{i=k}^m(-1)^{l-(i-k)}\bi{x+2k+i-k}{l-(i-k)}\bi{2k+2(i-k)}{i-k}
\\=&(-1)^m\sum_{i=k}^m(-1)^i\bi{m+n+i}{m-i}\bi{2i}{k+i}
\\=&\sum_{j=0}^l(-1)^j\bi{n+j-1}j\quad\ \t{(by [Su3, (3.2)])}
\\=&\sum_{j=0}^l\bi{-n}j=\sum_{j=0}^l\bi{l-x}j.\endalign$$
This concludes the proof. \qed
\medskip

\Remark\ 2.1. Lemma 2.1, an equivalent version of [Su3, (3.2)],
played a key role when the author established the following
curious identity in [Su3]:
$$\aligned&(x+m+1)\sum_{i=0}^m(-1)^i\bi{x+y+i}{m-i}\bi{y+2i}i
\\&=\sum_{i=0}^m\bi{x+i}{m-i}(-4)^i+(x-m)\bi xm,\endaligned\tag2.2$$
where $m$ is any nonnegative integer. The reader is referred to
[C], [CC], [EM], [MS] and [PP] for other proofs of
(2.2), and to [SW] for an extension of (2.2). In the case $x\in\{0,\ldots,l\}$ the right-hand side of
(2.1) turns out to be $2^{l-x}$, so (2.1) implies identity (3) in [C],
which has a nice combinatorial interpretation.

\medskip
\noindent{\it Proof of Theorem 1.1}. Let $n\in\Z$ and $n\gs2h$,
where $h=\lfloor(m-1)/2\rfloor$. Then
$n+1\gs m-1>m-2$. Suppose that (1.3) holds for all $k\in\Z$. Then,
for any given $k\in\Z$, we have
$$\align&\sum_{i=0}^{h}(-1)^i\bi{m-1-i}i\M{n+1-2i}{k-i}
\\=&\sum_{i=0}^{h}(-1)^i\bi{m-1-i}i
\(\M{n-2i}{k-i}+\M{n-2i}{k-1-i}\)
\\=&\sum_{i=0}^{h}(-1)^i\bi{m-1-i}i
\M{n-2i}{k-i}+\sum_{i=0}^{h}(-1)^i\bi{m-1-i}i\M{n-2i}{k-1-i}
\\=&2^{n-m+1}+\da_{m-2,\,n}\f{(-1)^k}2+\l(2^{n-m+1}+\da_{m-2,\,n}\f{(-1)^{k-1}}2\r)
=2^{(n+1)-m+1}.\endalign$$

In view of the above, it suffices to show (1.3) for
$n=2h$ and $k\in\{0,1,\ldots,m-1\}$. For any
$i\in\N$ with $i\ls h$, we have $k-i+m>n-2i$ since $n-m<0\ls k+i$, thus
$$\M{n-2i}{k-i}=\cases\bi{n-2i}{k-i}&\t{if}\ i\ls k,\\0&\t{if}\ i>k.\endcases$$
Let $x=m-1-n+k$, $y=n-2k$, and $\Sigma$ denote the left hand side
of (1.3). Then
$$\align\Sigma=&\sum_{i=0}^k(-1)^i\bi{m-1-i}i\bi{n-2i}{k-i}
\\=&\sum_{j=0}^k(-1)^{k-j}\bi{x+y+j}{k-j}\bi{y+2j}j
\\=&\sum_{j=0}^k\bi{k-x}j=\sum_{j=0}^k\bi{n-(m-1)}j
\endalign$$
with the help of Lemma 2.1. If $m$ is odd, then $n=m-1$ and hence
$\Sigma=\sum_{j=0}^k\bi0j=1=2^{n-m+1}$. If $m$ is even, then
$n=m-2$ and
$$\Sigma=\sum_{j=0}^k\bi{-1}j=\sum_{j=0}^k(-1)^j=\f{1+(-1)^k}2=2^{n-m+1}+\f{(-1)^k}2.$$
So we do have $\Sigma=2^{n-m+1}+\da_{m-2,\,n}(-1)^k/2$ as
required. \qed

\medskip
\noindent{\it Proof of Corollary 1.1}. Let $n\in\N$ and $n\gs2\lfloor (m-1)/2\rfloor$. By Theorem 1.1,
$$\sum_{i=0}^{\lfloor(m-1)/2\rfloor}(-1)^i\bi{m-1-i}i\M{n-2i}{\lfloor\f{k+n}2\rfloor-i}
=2^{n-m+1}+\da_{m-2,\,n}\f{(-1)^{\lfloor (k+n)/2\rfloor}}2.$$
If $m-2=n$, then $2\mid m$ and $(k+n)/2=(k+m)/2-1$. So (1.5) holds.

For $0\ls i\ls2\lfloor(m-1)/2\rfloor$, if $n-2i=0$ and $2\mid m$, then
we must have $n/2=i=\lfloor(m-1)/2\rfloor=m/2-1$.
Note also that
$$\sum_{i=0}^{\lfloor(m-1)/2\rfloor}(-1)^i\bi{m-1-i}i2^{m-1-2i}=m$$
by (1.60) of [G] or (4) of [C].
Therefore
$$\align&\sum_{i=0}^{\lfloor(m-1)/2\rfloor}(-1)^i\bi{m-1-i}iv_{n-2i}
\\=&m\sum_{i=0}^{\lfloor(m-1)/2\rfloor}(-1)^i\bi{m-1-i}iu_{n-2i}
-\sum_{i=0}^{\lfloor(m-1)/2\rfloor}(-1)^i\bi{m-1-i}i2^{n-2i}
\\&-\sum_{i=0}^{\lfloor(m-1)/2\rfloor}(-1)^i\bi{m-1-i}i
\da_{n-2i,\,0}\da_{(-1)^m,1}(-1)^{\lfloor k/2\rfloor}
\\=&-m\da_{m-2,\,n}\f{(-1)^{\lfloor (k+m)/2\rfloor}}2
-\da_{m-2,\,n}(-1)^{\lfloor (m-1)/2\rfloor}\f m2(-1)^{\lfloor k/2\rfloor}=0.
\endalign$$
This concludes the proof. \qed

\medskip
\noindent{\it Proof of Theorem 1.2}. i) Clearly $c_m(0)=1$.
As $\lfloor m/2\rfloor +\lfloor (m+1)/2\rfloor =m$, whether $m=1$ or not, we have
$$c_m\l(\l\lfloor  \f{m+1}2\r\rfloor \r)
=4\bi{\lfloor \f m2\rfloor }{\lfloor \f{m-1}2\rfloor}=4\bi{m-\lfloor \f{m+1}2\rfloor}
{\lfloor \f{m+1}2\rfloor -1}.$$
If $0<i<m/2$ then
$$\align c_m(i)=&\f{(m-i)!}{i!(m-2i)!}\cdot\f{m^2+m-2i}{(m-i)(m+1-2i)}
=\f{(m-i)!}{i!(m-2i)!}\l(\f{m-2i}{m-i}+\f{4i}{m+1-2i}\r)
\\=&\f{(m-1-i)!}{i!(m-1-2i)!}+4\f{(m-i)!}{(i-1)!(m+1-2i)!}=\bi{m-1-i}i+4\bi{m-i}{i-1}.
\endalign$$

Let $n\in\N$ and $n\gs2\lfloor (m+1)/2\rfloor$. Set $h=\lfloor(m-1)/2\rfloor$.
As $n>n-2\gs 2h$,
by Theorem 1.1 we have
$$\sum^{h}_{i=0}(-1)^i\bi{m-1-i}i\M{n-2i}
{k-i}=2^{n-m+1}$$
and
$$\sum^{h}_{i=0}(-1)^i\bi{m-1-i}i\M{n-2-2i}
{k-1-i}=2^{n-2-m+1}+\da_{m,n}\f{(-1)^{k-1}}2.$$
Therefore
$$\align0=&2^{n-m+1}-4\cdot2^{n-2-m+1}
\\=&\sum^{h}_{i=0}(-1)^i\bi{m-1-i}i\M{n-2i}{k-i}
\\&-4\(\sum^{h}_{i=0}(-1)^i\bi{m-1-i}i\M{n-2-2i}{k-1-i}
+\da_{m,n}\f{(-1)^{k}}2\)
\endalign$$
and hence
$$\align  2(-1)^{k}\da_{m,n}=&\M{n}{k}+\sum_{0<i<m/2}(-1)^i\bi{m-1-i}i\M{n-2i}{k-i}
\\&+4\sum^{h+1}_{j=1}(-1)^j\bi{m-j}{j-1}\M{n-2j}{k-j}
\\=&\M{n}{k}+\sum_{0<i<m/2}(-1)^i\(\bi{m-1-i}i+4\bi{m-i}{i-1}\)
\M{n-2i}{k-i}
\\&+(-1)^{\lfloor \f{m+1}2\rfloor }4\bi{m-\lfloor \f{m+1}2\rfloor }{\lfloor \f{m+1}2\rfloor-1}
\M{n-2\lfloor \f{m+1}2\rfloor }{k-\lfloor \f{m+1}2\rfloor}
\\=&\sum_{i=0}^{h+1}(-1)^ic_m(i)\M{n-2i}{k-i}.
\endalign$$
This proves the first part of Theorem 1.2.

ii) Observe that
$$\f m{m-i}\bi{m-i}i=2\bi{m-i}i-\bi{m-1-i}i\in\Z
\quad \t{for}\ i=0,\ldots,m-1.$$
In view of (1.2), it suffices to verify (1.8) in the case $n=2\lfloor m/2\rfloor$
and $0\ls k<m$. For any $i\in\N$ with $i\ls n/2=\lfloor m/2\rfloor$,
we have $k-i+m>n-2i$ if and only if $i=k=0$ and $m=n$, and thus
$$\m{n-2i}{k-i}=\cases0&\t{if}\ i>k,\ \t{or}\ i=k=0\ \&\ m=n,
\\\bi{n-2i}{k-i}&\t{otherwise}.\endcases$$
Therefore
$$\align&\sum_{i=0}^{\lfloor m/2\rfloor}(-1)^id_m(i)\m{n-2i}{k-i}
\\=&\sum_{i=0}^k(-1)^i\f{m}{m-i}\bi{m-i}i\bi{n-2i}{k-i}-\da_{k,0}\da_{m,n}
\\=&2\sum_{i=0}^k(-1)^i\bi{m-i}i\bi{n-2i}{k-i}
-\sum_{i=0}^k(-1)^i\bi{m-1-i}i\bi{n-2i}{k-i}-\da_{k,0}\da_{m,n}
\\=&2\sum_{j=0}^k\bi{n-m}j-\sum_{j=0}^k\bi{n-(m-1)}j-\da_{k,0}\da_{m,n}=\da_{m-1,\,n}(-1)^k
\endalign$$
with the help of Lemma 2.1.

The proof of Theorem 1.2 is now complete.
\qed

\medskip
\noindent{\it Proof of Corollary 1.2}. The case $p=2$ can be verified directly, so let $p>2$.
Clearly, (1.9) holds if and only if $p^{\lfloor(n-1)/(p-1)\rfloor}\mid\mm nrp$.
If $n\gs p$, then
$\mm nrp=\sum_{i=1}^{\lfloor p/2\rfloor}(-1)^{i-1}d_p(i)\mm{n-2i}{r-i}p$ by Theorem 1.2.
Since $p\mid d_p(i)$ for $i=1,\ldots,\lfloor p/2\rfloor$, we have the desired result
by induction on $n$. \qed

\heading{3. Proof of Theorem 1.3}\endheading

Let $m\in\Z^+$, $n\in\N$ and $r\in\Z$. Set
$$\bm{n}{r}{m}=\cases\M nr&\t{if}\ 2\mid m,
\\\m{n}{r}&\t{if}\ 2\nmid m;\endcases
\ \ \t{and}\ \ \bmb{n}{r}{m}=\cases\m nr&\t{if}\ 2\mid m,
\\\M{n}{r}&\t{if}\ 2\nmid m.\endcases
\tag3.1$$
Clearly
$$(-1)^r\bm nrm=\sum^n\Sb k=0\\m\mid k-r\endSb\bi nk(-1)^k,
\ \ (-1)^r\bmb nrm=\sum^n\Sb k=0\\m\mid k-r\endSb\bi
nk(-1)^{k+(k-r)/m}$$
and
$$\bm{n}rm+\bmb{n}rm=\M{n}{r}+\m{n}{r}=2\Mm nr{2m}
=2\bm nr{2m}.$$
Since $\M n{n-r}=\M nr$ and $\m n{n-r}=\m nr$, we have $\bm n{n-r}m=\bm nrm$
and also $\bmb n{n-r}m=\bmb nrm$.

\proclaim{Lemma 3.1} Let $l\in\N$, $m,n\in\Z^+$, $q\in\Z$ and
$(q,m)=1$. Then
$$U^{(q)}_l(m,n)=\sum_{r=0}^n\f{n-r}{1+\da_{r,0}}
\((-1)^{qr}\bm{2l}{l+qr}m-\f{\da_{l,0}}m\),\tag3.2$$
and $U^{(q)}_l(m,n)\in\Z$ if $1\ls l\ls\lfloor(m+1)/2\rfloor$.
When $(q,2m)=1$, we have
$$V^{(q)}_l(m,n)=\sum_{r=0}^n\f{n-r}{1+\da_{r,0}}
(-1)^{qr}\bmb{2l}{l+qr}m,\tag3.3$$
and also $V^{(q)}_l(m,n)\in\Z$ provided $1\ls l\ls\lfloor(m+1)/2\rfloor$.
\endproclaim
\Proof. Let $\rho=1$, or $\rho=-1$ and $(q,2m)=1$.
With the help of the identity
$$\f{2-x^n-x^{-n}}{2-x-x^{-1}}=n+\sum_{r=1}^n(n-r)(x^r+x^{-r})=\sum_{r=-n}^n(n-|r|)x^r$$
(cf. [GS, (2.2)], we have
$$\align&\sum\Sb \gamma^m=\rho\\\gamma\not=1\endSb
\f{2-\gamma^{qn}-\gamma^{-qn}}{2-\gamma^q-\gamma^{-q}}(2-\gamma-\gamma^{-1})^l
\\=&\sum\Sb \gamma^m=\rho\\\gamma\not=1\endSb\sum_{r=-n}^n(n-|r|)\gamma^{qr}
(2-\gamma-\gamma^{-1})^l
\\=&\sum_{r=-n}^n(n-|r|)\(\sum_{\gamma^m=\rho}\gamma^{qr}
(2-\gamma-\gamma^{-1})^l-\da_{\rho,1}\da_{l,0}\)
\\=&\sum_{r=-n}^n(n-|r|)\times
\cases(-1)^{qr}m\bm{2l}{l+qr}m-\da_{l,0}&\t{if}
\ \rho=1,\\(-1)^{qr}m\bmb{2l}{l+qr}m&\t{if}\ \rho=-1,\endcases
\endalign$$
where in the last step we apply Remark 2.1 of [Su2].
Note that $\bm{2l}{l+qr}m=\bm{2l}{l-qr}m$ and
$\bmb{2l}{l+qr}m=\bmb{2l}{l-qr}m$. So we have (3.2), also (3.3) holds if $(q,2m)=1$.

Suppose that $1\ls l\ls \lfloor (m+1)/2\rfloor$. In the case $m=1$, both
 $\bm{2l}l m=0$ and $\bmb{2l}l m=4$ are even. If $m>1$, then $l+m>2l$ and hence
 $\bm{2l}l m=\bmb{2l}l m=\bi{2l}l=2\bi{2l-1}l$.
Therefore $U_l^{(q)}(m,n)\in\Z$, and also $V_l^{(q)}(m,n)\in\Z$
when $(q,2m)=1$.

The proof of Lemma 3.1 is now complete.
\qed

\Remark\ 3.1. Let $q\in\Z$, $m\in\Z^+$ and $(q,m)=1$.
In view of (3.2), we have
$$U^{(q)}_0(m,n)=\f n2\l(1-\f1m\r)=\f{n(m-n)}{2m}\ \ \t{for}\ n=1,\ldots,m.\tag3.4$$
If $m>1$, then
$$U^{(q)}_l(m,1)=\f12\bi{2l}l=\bi{2l-1}l
\ \ \t{for}\ l=1,\cs,\l\lfloor  \f{m+1}2\r\rfloor .\tag3.5$$
When $(q,2m)=1$,  $V^{(q)}_0(m,n)=2U^{(q)}_0(2m,n)-U^{(q)}_0(m,n)=n/2$
for $n=1,\ldots,m$, and also
$V^{(q)}_l(m,1)=2U^{(q)}_l(2m,1)-U^{(q)}_l(m,1)=\bi{2l-1}l$
if $m>1$ and $1\ls l\ls\lfloor (m+1)/2\rfloor$.
\medskip

For positive integers $m$ and $n$, it is known that
$$\sum_{r=0}^{m-1}B_n\l(\f{x+r}m\r)=m^{1-n}B_n(x)\tag3.6$$
(due to Raabe), and
$$E_{n-1}(x)=\f2{n}\l(B_{n}(x)-2^{n}B_{n}\l(\f x2\r)\r).$$

\proclaim{Lemma 3.2} Let $n$ be a positive integer, and let $x$ be a real number. Then
$$nE_{n-1}(\{x\})=2(-1)^{\lfloor x\rfloor }
\l(B_n(\{x\})-2^nB_n\l(\l\{\f x2\r\}\r)\r).$$
\endproclaim
\Proof. Clearly $2\{x/2\}-\{x\}=\lfloor x\rfloor-2\lfloor x/2\rfloor\in\{0,1\}$.
If $2\mid\lfloor x\rfloor $, then
$$B_n(\{x\})-2^nB_n\l(\l\{\f x2\r\}\r)
=B_n(\{x\})-2^nB_n\l(\f{\{x\}}2\r)=\f n2E_{n-1}(\{x\}).
$$
By Raabe's formula (3.6),
$$B_n\l(\f{\{x\}}2\r)+B_n\l(\f{\{x\}+1}2\r)=2^{1-n}B_n(\{x\}).$$
So, if $2\nmid \lfloor x\rfloor$ then
$$\align B_n(\{x\})-2^nB_n\l(\l\{\f x2\r\}\r)
=&B_n(\{x\})-2^nB_n\l(\f{\{x\}+1}2\r)
\\=&B_n(\{x\})-2^n\l(2^{1-n}B_n(\{x\})-B_n\l(\f{\{x\}}2\r)\r)
\\=&2^nB_n\l(\f{\{x\}}2\r)-B_n(\{x\})
=-\f n2E_{n-1}(\{x\}).
\endalign$$
This concludes the proof. \qed

From Lemma 3.2 we have
\proclaim{Lemma 3.3} Let $p$ be an odd prime,  and let $m,n\in\Z^+$ and $p\nmid m$. Then
$$\aligned&\f{(-1)^{\lfloor pn/m\rfloor }}2E_{p-2}\l(\l\{\f{pn}m\r\}\r)+\f{2^{p-1}-1}p
\\\eq&B_{p-1}\l(\l\{\f{pn}{2m}\r\}\r)-B_{p-1}
-\(B_{p-1}\l(\l\{\f{pn}m\r\}\r)-B_{p-1}\)\ (\mo\ p).
\endaligned\tag3.7$$
\endproclaim
\Proof. By Lemma 3.2,
$$\align&(-1)^{\lfloor pn/m\rfloor }\f{p-1}2E_{p-2}\l(\l\{\f{pn}m\r\}\r)
+(2^{p-1}-1)B_{p-1}
\\=&B_{p-1}\l(\l\{\f{pn}m\r\}\r)-B_{p-1}
-2^{p-1}\(B_{p-1}\l(\l\{\f{pn}{2m}\r\}\r)-B_{p-1}\)\endalign$$
As $2^{p-1}\eq1\ (\mo\ p)$ by Fermat's little theorem,
 and $pB_{p-1}\eq-1\ (\mo\ p)$ by [IR, p.\,233], the desired (3.7) follows at once. \qed

\Remark\ 3.2. Let $p$ be an odd prime not dividing $m\in\Z^+$. By
[GS, pp.\,125--126] or [Su5, Corollary 2.1],
$$B_{p-1}\l(\l\{\f{pn}m\r\}\r)-B_{p-1}
\eq-\sum_{k=1}^{\lfloor pn/m\rfloor}\f1k\ (\mo\ p)\ \ \ \t{for}\ n=0,\ldots,m-1.$$
Combining this with (3.7) we get that
$$\align&\f{(-1)^{\lfloor pn/m\rfloor }}2E_{p-2}\l(\l\{\f{pn}m\r\}\r)+\f{2^{p-1}-1}p
\\\eq&\sum_{k=1}^{\lfloor pn/m\rfloor}\f1k-\sum_{k=1}^{\lfloor pn/(2m)\rfloor}\f1k
=\sum_{k=1}^{\lfloor pn/m\rfloor}\f{(-1)^{k-1}}k\ (\mo\ p)
\endalign$$
for every $n=0,\ldots,m-1$.
In light of Lemma 3.3, we can also deduce from (3) and (4) of [GS] the following congruences
with $n\in\Z^+$ and $(m,n)=1$.
$$\aligned
(-1)^{\lfloor pn/m\rfloor }E_{p-2}\l(\l\{\f{pn}m\r\}\r)
\eq\cases(\f 2n)\f4pP_{p-(\f 2p)}\ (\mo\ p)&\t{if}\ m=4,
\\(\f n5)\f5{p}F_{p-(\f 5p)}+\f{2^p-2}p\ (\mo\ p) &\t{if}\ m=5,
\\(\f 3{pn})\f6{p}S_{p-(\f 3p)}\ (\mo\ p)&\t{if}\ m=6,
\endcases\endaligned\tag3.8$$
where $(-)$ denotes the Jacobi symbol, and the sequences
$(F_k)_{k\in\N}$, $(P_k)_{k\in\N}$ and $(S_k)_{k\in\N}$ are
defined as follows:
$$\align&F_0=0,\ F_1=1,\ \t{and}\ F_{k+2}=F_{k+1}+F_k\ \t{for}\ k\in\N;
\\&P_0=0,\ P_1=1,\ \t{and}\ P_{k+2}=2P_{k+1}+P_k\ \t{for}\ k\in\N;
\\&S_0=0,\ S_1=1,\ \t{and}\ S_{k+2}=4S_{k+1}-S_k\ \t{for}\ k\in\N.
\endalign$$

\medskip
\noindent{\it Proof of Theorem 1.3}. Let $k\in\Z$.
By Theorem 1.2, for any integer $l\gs\lfloor m/2\rfloor$
we have
$$\align&\sum_{i=0}^{\lfloor m/2\rfloor}(-1)^ia_m(i)\((-1)^k\bm{2l-2i}{k+l-i}m
-\f{\da_{l-i,\,0}}m\)
\\=&(1+\da_{(-1)^m,1})(-1)^l\da_{l,\lfloor m/2\rfloor}
-\f{\da_{l,\lfloor m/2\rfloor}}m(-1)^la_m\l(\l\lfloor \f m2\r\rfloor \r)=0;
\endalign$$
also $$\sum_{j=0}^{\lfloor (m+1)/2\rfloor}(-1)^jb_m(j)\bmb{2l-2j}{k+l-j}m=0$$
for all integers $l\gs\l\lfloor (m+1)/2\r\rfloor$.
This, together with Lemma 3.1, yields (1.13), and also (1.14) in the case $(q,2m)=1$.

By Lemma 3.1, $U^{(q)}_l(m,n)\in\Z$ for every $l=1,\ldots,\lfloor(m+1)/2\rfloor$;
by Theorem 1.2, $a_m(i)\in\Z$ if $0<i\ls\lfloor m/2\rfloor$.
Thus, in view of (1.13), we have $U^{(q)}_l(m,n)\in\Z$ for
each $l=1,2,3,\ldots$. If $(q,2m)=1$, then
$V^{(q)}_l(m,n)=2U^{(q)}_l(2m,n)-U^{(q)}_l(m,n)\in\Z$ for all $l\in\Z^+$.

Now assume that $m\nmid n$, and let $p$ be an odd prime with $p\eq\pm q\ (\mo\ 2m)$.
By Lemma 3.3 and (1.12),
$$\align&(-1)^{\lfloor pn/m\rfloor}E_{p-2}\l(\l\{\f{pn}m\r\}\r)+\f{2^p-2}p
\\\eq&2\(B_{p-1}\l(\l\{\f{pn}{2m}\r\}\r)-B_{p-1}\)
-2\(B_{p-1}\l(\l\{\f{pn}{m}\r\}\r)-B_{p-1}\)
\\\eq&\f{2m}p\l(U^{(q)}_p(2m,n)-1\r)-\f mp\l(U^{(q)}_p(m,n)-1\r)
=\f mp\l(V^{(q)}_p(m,n)-1\r)\ (\mo\ p).
\endalign$$
This proves (1.16). We are done. \qed
\medskip

We can also prove (1.13) by determining the characteristic polynomial
$$f_m(x):=\prod_{0<k\ls \lfloor m/2\rfloor}\l(x-\l(2-e^{2\pi ik/m}-e^{-2\pi ik/m}\r)\r)\tag3.9$$
of the recurrence $(U_l^{(q)}(m,n))_{l\in\N}$ of order $\lfloor m/2\rfloor$. If $m$ is even, then
$$\align f_m(x)=&\prod_{0<k\ls m/2}\l(x-2-e^{2\pi i(m/2-k)/m}-e^{-2\pi i(m/2-k)/m}\r)
\\=&\prod_{j=0}^{m/2-1}\l(x-2-2\cos\f{2j\pi}m\r)
=(x-4)\prod\Sb 0<k<m\\2\mid k\endSb\l(x-4\cos^2\f{k\pi}{2m}\r).
\endalign$$
If $m$ is odd, then
$$\align f_m(x)=&\prod_{0<j\ls(m-1)/2}\l(x-2-e^{2\pi i(m-2j)/(2m)}-e^{-2\pi i(m-2j)/(2m)}\r)
\\=&\prod\Sb 0<k<m\\2\nmid k\endSb\l(x-2-2\cos\f{2k\pi}{2m}\r)
=\prod\Sb 0<k<m\\2\nmid k\endSb\l(x-4\cos^2\f{k\pi}{2m}\r).
\endalign$$
So $f_m(x)$ can be determined with the help of the following lemma.

\proclaim{Lemma 3.4} Let $n$ be any positive integer. Then
$$\prod\Sb 0<k<n\\2\mid k-\da\endSb\l(x-4\cos^2\f{k\pi}{2n}\r)=\cases C_n(x)&\t{if}\ \da=0,\\D_n(x)&\t{if}\ \da=1,
\endcases\tag3.10$$
where
$$C_n(x)=\sum_{i=0}^{\lfloor(n-1)/2\rfloor}(-1)^i\bi{n-1-i}ix^{\lfloor(n-1)/2\rfloor-i}\tag3.11$$
and
$$D_n(x)=\sum_{i=0}^{\lfloor n/2\rfloor}(-1)^i\f n{n-i}\bi{n-i}ix^{\lfloor n/2\rfloor-i}.\tag3.12$$
\endproclaim
\Proof. It is well known that $\cos(n\theta)=T_n(\cos\theta)$ and
$\sin(n\theta)=\sin\theta\cdot U_{n-1}(\cos\theta)$, where the Chebyshev polynomials
$T_n(x)$ and $U_{n-1}(x)$ are given by
$$T_n(x)=\f n2\sum_{i=0}^{\lfloor n/2\rfloor}(-1)^i\f{(n-1-i)!}{i!(n-2i)!}(2x)^{n-2i}$$
and
$$U_{n-1}(x)=\sum_{i=0}^{\lfloor (n-1)/2\rfloor}(-1)^i\f{(n-1-i)!}{i!(n-1-2i)!}(2x)^{n-1-2i}.$$
If $n$ is even, then $T_n(x)=D_n(4x^2)/2$ and $U_{n-1}(x)=2xC_n(4x^2)$;
if $n$ is odd, then $T_n(x)=xD_n(4x^2)$ and $U_{n-1}(x)=C_n(4x^2)$.

As $U_{n-1}(\cos\f{k\pi}{2n})=0$ for those even $0<k<n$, the $2\lfloor(n-1)/2\rfloor$ distinct numbers
$\pm\cos\f{k\pi}{2n}\ (0<k<n,\,2\mid k)$ are zeroes of the polynomial $C_n(4x^2)$ of degree $2\lfloor(n-1)/2\rfloor$.
Similarly, since $T_{n}(\cos\f{k\pi}{2n})=0$ for those odd $0<k<n$, the $2\lfloor n/2\rfloor$ distinct numbers
$\pm\cos\f{k\pi}{2n}\ (0<k<n,\,2\nmid k)$ are zeroes of the polynomial $D_n(4x^2)$ of degree $2\lfloor n/2\rfloor$.
So
$$C_n(4x^2)=\prod\Sb 0<k<n\\2\mid k\endSb\(4x^2-4\cos^2\f{k\pi}{2n}\)
\ \t{and}\ D_n(4x^2)=\prod\Sb 0<k<n\\2\nmid k\endSb\(4x^2-4\cos^2\f{k\pi}{2n}\).$$
Therefore (3.10) holds. \qed

\Remark\ 3.3. For each $n\in\Z^+$, by Lemma 3.4 we have
$$C_n(x)=\prod_{0<k<n/2}\l(x-4\cos^2\f{k\pi}n\r)=\prod_{d\mid n}A_d(x),$$
where
$$A_d(x)=\prod\Sb 0<c<d/2\\(c,d)=1\endSb\l(x-4\cos^2\f{c\pi}d\r).\tag3.13$$
Applying the M\"obius inversion formula we obtain
$$A_n(x)=\prod_{d\mid n}C_d(x)^{\mu(n/d)},\tag3.14$$
which makes the polynomial $A_n(x)$ (introduced in [Su2]) computable.

\heading{4. Appendix: Explicit Values of $a_m(i)$ and $b_m(j)$ for $2\ls m\ls 12$}\endheading

\medskip
\centerline{Table 1: Values of $a_m(i)$ with $2\ls m\ls 12$}
\smallskip
\centerline{\vbox{\offinterlineskip
\halign{\vrule#&\ \ #\ \hfill   &&\vrule#&\ \ \hfill#\ \
\cr\nh\cr \hh\om\om\om\om\om\om\om
\cr &$\st {\ \ i}{m\ \ }$& &$1$& &$2$& &$3$& &$4$& &$5$& &$6$&
\cr \hh\om\om\om\om\om\om\om
\cr\nh\cr \hh\om\om\om\om\om\om\om
\cr & $2$ && $4$ &&$ $ && $ $ && $ $ && $ $ && $ $&
\cr \hh\om\om\om\om\om\om\om
\cr\nh\cr \hh\om\om\om\om\om\om\om
\cr & $3$ && $3$ && $ $ && $ $ && $ $ && $ $ && $ $&
\cr \hh\om\om\om\om\om\om\om
\cr \nh\cr \hh\om\om\om\om\om\om\om
\cr & $4$ && $6$ && $8$ && $ $ && $ $ && $ $ && $ $&
\cr \hh\om\om\om\om\om\om\om
\cr \nh\cr \hh\om\om\om\om\om\om\om
\cr & $5$ && $5$ && $5$ && $ $ && $ $ && $ $ && $ $&
\cr \hh\om\om\om\om\om\om\om
\cr \nh\cr \hh\om\om\om\om\om\om\om
\cr & $6$ && $8$ && $19$ && $12$ && $ $ && $ $ && $ $&
\cr \hh\om\om\om\om\om\om\om
\cr \nh\cr \hh\om\om\om\om\om\om\om
\cr & $7$ && $7$ && $14$ && $7$ && $ $ && $ $ && $ $&
\cr \hh\om\om\om\om\om\om\om
\cr \nh\cr \hh\om\om\om\om\om\om\om
\cr & $8$ && $10$ && $34$ && $44$ && $16$ && $ $ && $ $&
\cr \hh\om\om\om\om\om\om\om
\cr \nh\cr\hh\om\om\om\om\om\om\om
\cr & $9$ && $9$ && $27$ && $30$ && $9$ && $ $ && $ $&
\cr \hh\om\om\om\om\om\om\om
\cr \nh\cr \hh\om\om\om\om\om\om\om
\cr & $10$ && $12$ && $53$ && $104$ && $85$ && $20$ && $ $&
\cr \hh\om\om\om\om\om\om\om
\cr \nh\cr \hh\om\om\om\om\om\om\om
\cr & $11$ && $11$ && $44$ && $77$ && $55$ && $11$ && $ $&
\cr \hh\om\om\om\om\om\om\om
\cr \nh\cr \hh\om\om\om\om\om\om\om
\cr & $12$ && $14$ && $76$ && $200$ && $259$ && $146 $ && $24$&
\cr \hh\om\om\om\om\om\om\om
\cr \nh\cr }}}
\smallskip
\bigskip
\bigskip
\centerline{Table 2: Values of $b_m(j)$ with $2\ls m\ls 12$}
\smallskip
\centerline{\vbox{\offinterlineskip
\halign{\vrule#&\ \ #\ \hfill   &&\vrule#&\ \ \hfill#\ \
\cr\nh\cr \hh\om\om\om\om\om\om\om
\cr &$\st {\ \ j}{m\ \ }$& &$1$& &$2$& &$3$& &$4$& &$5$& &$6$&
\cr \hh\om\om\om\om\om\om\om
\cr\nh\cr \hh\om\om\om\om\om\om\om
\cr & $2$ && $2$ &&$ $ && $ $ && $ $ && $ $ && $ $&
\cr \hh\om\om\om\om\om\om\om
\cr\nh\cr \hh\om\om\om\om\om\om\om
\cr & $3$ && $5$ && $4$ && $ $ && $ $ && $ $ && $ $&
\cr \hh\om\om\om\om\om\om\om
\cr \nh\cr \hh\om\om\om\om\om\om\om
\cr & $4$ && $4$ && $2$ && $ $ && $ $ && $ $ && $ $&
\cr \hh\om\om\om\om\om\om\om
\cr \nh\cr \hh\om\om\om\om\om\om\om
\cr & $5$ && $7$ && $13$ && $4$ && $ $ && $ $ && $ $&
\cr \hh\om\om\om\om\om\om\om
\cr \nh\cr \hh\om\om\om\om\om\om\om
\cr & $6$ && $6$ && $9$ && $2$ && $ $ && $ $ && $ $&
\cr \hh\om\om\om\om\om\om\om
\cr \nh\cr \hh\om\om\om\om\om\om\om
\cr & $7$ && $9$ && $26$ && $25$ && $4$ && $ $ && $ $&
\cr \hh\om\om\om\om\om\om\om
\cr \nh\cr \hh\om\om\om\om\om\om\om
\cr & $8$ && $8$ && $20$ && $16$ && $2$ && $ $ && $ $&
\cr \hh\om\om\om\om\om\om\om
\cr \nh\cr\hh\om\om\om\om\om\om\om
\cr & $9$ && $11$ && $43$ && $70$ && $41$ && $4$ && $ $&
\cr \hh\om\om\om\om\om\om\om
\cr \nh\cr \hh\om\om\om\om\om\om\om
\cr & $10$ && $10$ && $35$ && $50$ && $25$ && $2$ && $ $&
\cr \hh\om\om\om\om\om\om\om
\cr \nh\cr \hh\om\om\om\om\om\om\om
\cr & $11$ && $13$ && $64$ && $147$ && $155$ && $61$ && $4$&
\cr \hh\om\om\om\om\om\om\om
\cr \nh\cr \hh\om\om\om\om\om\om\om
\cr & $12$ && $12$ && $54$ && $112$ && $105$ && $36 $ && $2$&
\cr \hh\om\om\om\om\om\om\om
\cr \nh \cr}}}

\bigskip
\Ack. The author thanks the referee for his/her helpful comments.

\enddocument